# Further Equivalences and Semiglobal Versions of Integral Input to State Stability

David Angeli
Dip. Sistemi e Informatica
University of Florence
50139 Firenze, Italy
angeli@dsi.unifi.it

Eduardo Sontag[*]
Department of Mathematics
Rutgers University
New Brunswick, NJ 08903, USA
sontag@hilbert.rutgers.edu

Yuan Wang[†]
Department of Mathematical Sciences
Florida Atlantic University
Boca Raton, FL 33431, USA
ywang@control.math.fau.edu

**Abstract**

This paper continues the study of the integral input-to-state stability (iISS) property. It is shown that the iISS property is equivalent to one which arises from the consideration of mixed norms on states and inputs, as well as to the superposition of a "bounded energy bounded state" requirement and the global asymptotic stability of the unforced system. A semiglobal version of iISS is shown to imply the global version, though a counterexample shows that the analogous fact fails for input to state stability (iss). The results in this note complete the basic theoretical picture regarding iISS and ISS.

## 1 Introduction and Basic Definitions

We consider continuous time nonlinear systems of the following form:

$$\dot{x} = f(x, u) \tag{1}$$

with states $x$ evolving in $\mathbb{R}^n$ and inputs $u$ taking values in $\mathbb{R}^m$. Inputs $u$ are measurable, locally essentially bounded functions of time. The map $f : \mathbb{R}^n \times \mathbb{R}^m \to \mathbb{R}^n$ is assumed to satisfy $f(0,0) = 0$, and to be locally Lipschitz. Given any state $\xi \in \mathbb{R}^n$ and any input $u : [0, \infty) \to \mathbb{R}^m$ we denote by $x(t, \xi, u)$ the unique maximal solution of the system (1), which is defined on some maximal interval $[0, t_{\max}(\xi, u))$. The system is said to be *forward complete* if $t_{\max}(\xi, u) = +\infty$ for all $\xi$ and $u$. We use the notation $|\zeta|$ for Euclidean norm of vectors $\zeta$, and $\|u\|_\infty$ for (essential) supremum of a function of time.

This paper continues the study of the input-to-state stability (iss) property (cf. [4, 6, 7, 8, 9, 10, 11, 15, 14, 16, 17, 19, 20, 21, 22]) as well as its variant, the *integral* input-to-state stability (IISS) property (cf. [1, 2, 12, 18]). Recall that the first of these is the natural extension to nonlinear systems (under arbitrary coordinate changes in states and inputs) of the notion of external stability known as "finiteness of $\mathcal{L}^1$ gain" (that is, finite operator norm from $\mathcal{L}^\infty$ to $\mathcal{L}^\infty$), while the second one extends to nonlinear systems the notion of finite "$\mathcal{H}_2$ gain" (that is, finite operator norm from $\mathcal{L}^2$ to $\mathcal{L}^\infty$).

---

[*]Supported in part by US Air Force Grant F49620-98-1-0242
[†]Supported in part by NSF Grant DMS-9457826

It was shown in [18] that the ISS property is equivalent to the following "integral to integral" stability concept: the system (1) is forward complete, and there are some $\alpha$, $\chi$, and $\sigma \in \mathcal{K}_\infty$ such that, for all initial states $\xi$ and inputs $u$, the following estimate holds for all $t \in [0, t_{\max}(\xi, u))$:

$$\int_0^t \alpha(|x(s, \xi, u)|)\, ds \;\leq\; \chi(|\xi|) \;+\; \int_0^t \sigma(|u(s)|)\, ds\,. \tag{2}$$

In particular, for example, a system which has "finite $\mathcal{H}_\infty$ norm" (i.e., operator norm from $\mathcal{L}^2$ to $\mathcal{L}^2$), meaning that an estimate such as

$$\left(\int_0^t (|x(s, \xi, u)|)^q\, ds\right)^{\frac{1}{q}} \;\leq\; \left(|\xi|^p + \int_0^t \sigma(|u(s)|)^p\, ds\right)^{\frac{1}{p}} \tag{3}$$

holds with $p = q = 2$, is necessarily ISS. (Actually, in a sense every ISS system satisfies an estimate of this form, under coordinate changes; see [5].) The estimate (3) with $p = q = 2$ (or, more generally, for any $p = q$) gives rise to an estimate as in (2) simply by raising both sides to the $p$th power. However, *mixed* norms $p \neq q$ give rise, after raising to the $p$th power, to the more general type of estimate:

$$\gamma\left(\int_0^t \alpha(|x(s, \xi, u)|)\, ds\right) \;\leq\; \chi(|\xi|) \;+\; \int_0^t \sigma(|u(s)|)\, ds \tag{4}$$

with all comparison functions of class $\mathcal{K}_\infty$. Alternative ways of stating such a condition arise by taking $\gamma^{-1}$, leading to an estimate such as

$$\int_0^t \alpha(|x(s, \xi, u)|)\, ds \;\leq\; \gamma\left(\chi(|\xi|) + \int_0^t \sigma(|u(s)|)\, ds\right) \tag{5}$$

(we wrote $\gamma^{-1}$ again as $\gamma$), which could also be written as

$$\int_0^t \alpha(|x(s, \xi, u)|)\, ds \;\leq\; \tilde{\chi}(|\xi|) \;+\; \tilde{\gamma}\left(\int_0^t \sigma(|u(s)|)\, ds\right) \tag{6}$$

(just let $\tilde{\chi}(r) := \gamma(2\chi(r))$ and $\tilde{\gamma}(r) := \gamma(2r)$), or even as:

$$\int_0^t \alpha(|x(s, \xi, u)|)\, ds \;\leq\; \chi\left(|\xi| + \int_0^t \sigma(|u(s)|)\, ds\right) \tag{7}$$

holding for some $\alpha$, $\chi$, and $\sigma \in \mathcal{K}_\infty$ (take $\chi := \tilde{\chi} + \tilde{\gamma}$). Note that an estimate of this type in turn implies again an estimate as in (4), if one takes $\gamma := \chi^{-1}$ (and "$\chi$" is the identity in (4)).

Given the apparent similarity between on the one hand (2) and on the other (4) and its equivalent versions (5)-(7), it seems natural to conjecture that this latter property is also equivalent to ISS. Surprisingly, however, and this is one of the main results of this paper, this equivalence turns out to be false, and (4) is in fact equivalent to *integral input-to-state stability* (IISS). Recall that a system is said to be IISS if there exist functions $\alpha, \sigma$ of class $\mathcal{K}_\infty$, and $\beta$ of class $\mathcal{KL}$, such that

$$\alpha(|x(t,\xi,u)|) \leq \beta(|\xi|,t) + \int_0^t \sigma(|u(s)|)\,ds \tag{8}$$

holds along all trajectories. Observe that any IISS system is necessarily forward complete, because $\alpha(|x(t,\xi,u)|)$ is bounded by $\beta(|\xi|,0) + \int_0^T \sigma(|u(s)|)\,ds$ on any interval $[0,T)$, with $T < t_{\max}(\xi,u)$, on which the trajectory is defined (so maximal trajectories stay bounded, and are therefore everywhere defined).

The concept of IISS is very natural; among other characterizations, it was shown in [2] that IISS is equivalent to the existence of a proper and positive definite smooth function $V$ which satisfies a dissipation inequality of the following kind:

$$DV(x)f(x,u) \leq -\rho(|x|) + \sigma(|u|) \quad \forall x \in \mathbb{R}^n,\ \forall u \in \mathbb{R}^m, \tag{9}$$

for some positive definite function $\rho$ and some $\sigma$ of class $\mathcal{K}_\infty$. (In contrast, the strictly stronger ISS property is equivalent to the existence of a dissipation inequality of this general form, but where $\rho$ is required to be class $\mathcal{K}_\infty$.)

We may weaken the requirements in the IISS definition by not asking that the effect of initial conditions decay, replacing $\beta(|\xi|,t)$ by just an upper bound $\gamma(|\xi|)$, and even by allowing an additional additive constant. This leads to the following notion: a system is *uniformly bounded energy bounded state* (UBEBS) if, for some $\alpha$, $\gamma$, and $\sigma \in \mathcal{K}_\infty$, and some positive constant $c$, the following estimate holds along all trajectories:

$$\alpha(|x(t,\xi,u)|) \leq \gamma(|\xi|) + \int_0^t \sigma(|u(s)|)\,ds + c. \tag{10}$$

Another apparently weaker notion of stability, which will be used in Section 3 in order to deal with semi-global versions of IISS, is given by estimates of this type, which mix integral and sup norms:

$$\alpha(|x(t,\xi,u)|) \leq \beta(|\xi|,t) + \int_0^t \sigma(|u(s)|)\,ds + \gamma(\|u_{[0,t]}\|_\infty), \tag{11}$$

understood as holding for some $\beta$ of class $\mathcal{KL}$ and some $\alpha$, $\sigma$, and $\gamma$ of class $\mathcal{K}_\infty$.

Our main equivalence results, to be proved in Section 2, can be summarized as follows. In addition to the "integral to integral stability" equivalences, we also state a sort of "separation" theorem for IISS, which allows to decompose the property into global asymptotic stability of the zero-input system $\dot{x} = f(x,0)$ (the 0-GAS property) plus the "bounded input-energy bounded state" property UBEBS.

**Theorem 1** *For any system (1), the following facts are equivalent:*

1. *The system is* IISS.

2. *The system is* 0-GAS *and* UBEBS.

3. *The system is forward complete and satisfies an estimate as in (7).*

4. *The system satisfies an estimate as in (11).*

We turn now to a "semiglobal" version of the IISS property, in which estimates are only required to hold for bounded initial states and inputs.

**Definition 1.1** A system (1) is *semiglobally* IISS if for each $M > 0$ there are functions $\beta_M \in \mathcal{KL}$, and $\gamma_M$ and $\alpha_M$ in $\mathcal{K}_\infty$, such that the following estimate:

$$\alpha_M(|x(t,\xi,u)|) \leq \beta_M(|\xi|,t) + \int_0^t \sigma_M(|u(s)|)\,ds \qquad (12)$$

holds for all initial states $\xi$ and inputs $u$ such that $|\xi| \leq M$ and $\|u\|_\infty \leq M$, for all $t \in [0, t_{\max}(\xi, u))$. □

Note that such a system is clearly O-GAS, because when $u \equiv 0$, $\beta_{|\xi|}(|\xi|,t) \to 0$ as $t \to \infty$ implies that $x(t,\xi,0) \to 0$ (attractivity), and applying with $M = 1$ one establishes stability.

The main result in that respect will be as follows.

**Theorem 2** *A system (1) is semiglobally IISS if and only if it is IISS.*

This is proved in Section 3. Interestingly, the respective result does not hold for the ISS property, in so far as initial states are concerned; see Section 4.

## 2 Proof of Theorem 1

We will first show that $1 \Leftrightarrow 2$. It is obvious that $1 \Rightarrow 2$, but the converse requires the following steps: first we establish the existence of a lower-semicontinuous Lyapunov-like function $V$, under the assumption that an UBEBS-like estimate holds, and then we combine this $V$ with a function as in the characterization of O-GAS given in [2] to obtain a non-smooth dissipation inequality; the final step is to show that the IISS property can be deduced from this inequality.

Next we establish that $1 \Rightarrow 3$. This implication makes essential use of the Lyapunov characterization of IISS given in [2].

We then turn to showing that $3 \Rightarrow 2$. The proof of the implication is heavily based upon the Lyapunov characterization of forward completeness which was recently given in [3]. (As a matter of fact, the original motivation for that paper was in trying to provide the main technical step required in this proof.) The fact that (7) together with forward completeness is a sufficient condition for O-GAS was already shown in [19] (see the proof of Theorem 1 in that paper, applied to the special case when $u = 0$). Hence, we are only left to show that a UBEBS estimate holds.

Since it is clearly true that 1 implies 4, we are only left to show that the converse is also true; we do this in Section 2.4.

### 2.1  $1 \Longleftrightarrow 2$

The implication $1 \Rightarrow 2$ follows easily simply considering that $\beta(|\xi|,t) \leq \beta(|\xi|,0)$ and recalling that $\beta(\cdot,0)$ is a $\mathcal{K}$ function. The converse implication is more interesting.

As a preliminary step, we show that, if a system is O-GAS and UBEBS, we may always reduce to the special case $c = 0$. Moreover, we show that a weaker estimate is also possible.

**Lemma 2.1** Suppose that a system (1) is O-GAS. Then the following properties are equivalent:

- The system satisfies along all trajectories an estimate of the following type, for suitable maps of class $\mathcal{K}_\infty$:
$$\alpha_1(|x(t,\xi,u)|) \leq \alpha_3(|\xi|) + \chi\left(\int_0^t \kappa(|u(s)|)\,ds\right). \tag{13}$$

- The system satisfies along all trajectories an UBEBS-like estimate with $c = 0$:
$$\alpha(|x(t,\xi,u)|) \leq \gamma(|\xi|) + \int_0^t \sigma_1(|u(s)|)\,ds. \tag{14}$$

- The system is UBEBS.

*Proof.* If an estimate of the type (13) holds along all trajectories, we simply introduce $\alpha_4(r) := \chi^{-1}(r/2)$, so that (14) holds with $\alpha = \alpha_4 \circ \alpha_1$ and $\gamma(r) = \alpha_4(2\alpha_3(r))$. Clearly, if an estimate of type (14) holds, then the system is UBEBS (just take $c = 0$). Thus, all that we need to prove is that any O-GAS and UBEBS system satisfies an estimate of type (13).

By virtue of Lemma 4.10 in [2], we have that O-GAS implies the existence of a smooth function $V : \mathbb{R}^n \to \mathbb{R}_{\geq 0}$, two class $\mathcal{K}_\infty$ functions $\alpha_i$ such that (18) holds, and some $\theta, \delta$ of class $\mathcal{K}_\infty$, so that
$$DV(x)\,f(x,u) \leq \theta(|x|)\,\delta(|u|) \quad \forall x \in \mathbb{R}^n,\ \forall u \in \mathbb{R}^m. \tag{15}$$
Taking the integral of this inequality in both sides yields the following estimate:
$$V(x(t,\xi,u)) \leq V(\xi) + \int_0^t \theta(|x(s,\xi,u)|)\,\delta(|u(s)|)\,ds \tag{16}$$
along all possible solutions corresponding to initial states $\xi$ and controls $u$. We now exploit the UBEBS property in order to show that we can always reduce ourselves to the case $c = 0$.

In general, for any class-$\mathcal{K}$ map $\alpha$, we have that $\alpha(a + b + c) \leq \alpha(3a) + \alpha(3b) + \alpha(3c)$. Applying this observation to (10), we obtain the following estimate along all solutions:
$$\theta(|x(s,\xi,u)|) \leq \tilde{\alpha}(3\gamma(|\xi|)) + \tilde{\alpha}\left(\int_0^s 3\sigma(|u(\tau)|)\,d\tau\right) + \tilde{\alpha}(3c), \tag{17}$$
where $\tilde{\alpha}(r) := \theta \circ \alpha^{-1}(3r)$. Majorizing $\theta(|x(s)|)$ in the right hand side of (16) according to (17), we obtain:
$$V(x(t,\xi,u)) - V(\xi)$$
$$\leq \int_0^t \left\{\tilde{\alpha}(\gamma(|\xi|)) + \tilde{\alpha}\left(\int_0^s \sigma(|u(\tau)|)\,d\tau\right) + \tilde{\alpha}(c)\right\}\delta(|u(s)|)\,ds$$
$$\leq \int_0^t \left\{\tilde{\alpha}(\gamma(|\xi|)) + \tilde{\alpha}\left(\int_0^t \sigma(|u(\tau)|)\,d\tau\right) + \tilde{\alpha}(c)\right\}\delta(|u(s)|)\,ds$$
$$= \left\{\tilde{\alpha}(\gamma(|\xi|)) + \tilde{\alpha}\left(\int_0^t \sigma(|u(s)|)\,ds\right) + \tilde{\alpha}(c)\right\}\int_0^t \delta(|u(s)|)\,ds$$
$$\leq \left[\tilde{\alpha}(\gamma(|\xi|))\right]^2 + \chi\left(\int_0^t \kappa(|u(s)|)\,ds\right)$$
where $\kappa(r) = \max\{\sigma(r), \delta(r)\}$ and $\chi(r) = r^2 + \tilde{\alpha}(c)r + \tilde{\alpha}(r)r$. Now let $\alpha_1$ and $\alpha_2$ be such that (18) holds. Then we have an estimate of the type (13) holding along all trajectories, with with $\alpha_3 := \alpha_2 + [\tilde{\alpha}(\gamma(\cdot))]^2$. ∎

Thus, we assume from now on that an estimate (14) holds. The next step is to obtain a Lyapunov-like property.

**Lemma 2.2** Suppose that system (1) satisfies an estimate (14). Then, there exist functions $\alpha_1$, $\alpha_2$, and $\sigma_1$ of class $\mathcal{K}_\infty$, and a lower semicontinuous function $V : \mathbb{R}^n \to \mathbb{R}_{\geq 0}$ such that

$$\alpha_1(|x|) \leq V(x) \leq \alpha_2(|x|) \tag{18}$$

holds for all $x$, so that along all trajectories the following estimate is satisfied:

$$V(x(t, \xi, u)) - V(\xi) \leq \int_0^t \sigma_1(|u(s)|) \, ds. \tag{19}$$

*Proof.* Take as $V$ the following function:

$$V(\xi) := \sup_{t \geq 0, u(\cdot)} \left\{ \alpha(|x(t, \xi, u)|) - \int_0^t \sigma_1(|u(s)|) \, ds \right\}. \tag{20}$$

Lower-semicontinuity of $V$ follows by a routine argument from the continuity of $\alpha(|x(t, \cdot, u)|)$. By definition of UBEBS in (10), $V$ is finite-valued and (18) is satisfied with $\alpha_1 = \alpha$ and $\alpha_2 = \gamma$.

We show next that, along trajectories of the system, $V$ satisfies (19). In fact

$$\begin{aligned}
V(x(t, \xi, u)) &= \sup_{\tau \geq 0, v(\cdot)} |x(\tau, x(t, \xi, u), v)| - \int_0^t \sigma_1(|v(s)|) \, ds \\
&= \sup_{\tau \geq 0, v(\cdot)} |x(\tau + t, \xi, u\#v)| - \int_0^t \sigma_1(|v(s)|) \, ds \\
&= \left[ \sup_{\tau \geq 0, v(\cdot)} |x(\tau + t, \xi, u\#v)| - \int_0^{t+\tau} \sigma_1(|u\#v(s)|) \, ds \right] + \int_0^t \sigma_1(|u(s)|) \, ds \\
&\leq \left[ \sup_{\tilde{\tau} \geq 0, \tilde{v}(\cdot)} |x(\tilde{\tau}, \xi, \tilde{v})| - \int_0^{\tilde{\tau}} \sigma_1(|\tilde{v}(s)|) \, ds \right] + \int_0^t \sigma_1(|u(s)|) \, ds \\
&\leq V(\xi) + \int_0^t \sigma_1(|u(s)|) \, ds,
\end{aligned}$$

where $u\#v$ denotes the concatenation input: $u\#v(t) = u(t)$ for $t < \tau$ and $v(t - \tau)$ otherwise. ∎

To complete the proof that the system is IISS, we take a function $V_1$ as in Lemma 2.2. In [2], it is shown that, for each O-GAS system there exists some smooth, positive definite, and "semi-proper" function $V_2 : \mathbb{R}^n \to \mathbb{R}_{\geq 0}$ so that, for some $\sigma_2 \in \mathcal{K}_\infty$ and some continuous positive definite function $\rho : \mathbb{R}_{\geq 0} \to \mathbb{R}_{\geq 0}$,

$$DV_2(x) f(x, u) \leq -\rho(|x|) + \sigma_2(|u|) \tag{21}$$

holds for all $x \in \mathbb{R}^n$ and all $u \in \mathbb{R}^m$. (Semiproper was defined in that paper as: for each $r$ in the range of $V_2$, the sublevel set $\{x | V_2(x) \leq r\}$ is compact.) Consider the function $V := V_1 + V_2$. This function is such that (18) holds for all $x$, for suitable $\alpha_1$, $\alpha_2$ of class $\mathcal{K}_\infty$, and also there is a $\sigma$ of class $\mathcal{K}_\infty$ (namely, we may pick $\sigma = \sigma_1 + \sigma_2$) and a continuous positive definite function $\rho$, such that

$$V(x(t, \xi, u)) \leq V(\xi) + \int_0^t \sigma(|u(s)|) \, ds - \int_0^t \rho(|x(s, \xi, u)|) \, ds \tag{22}$$

holds along all trajectories. (Note that $V(x(t, \xi, u))$ is Lebesgue measurable as a function of $t$, because $V$ is lower semicontinuous, so the integral makes sense.) Thus the proof will be completed once that we show the following non-smooth version of the sufficiency condition for iiss established in [2].

**Lemma 2.3** Consider a system (1), and suppose that there exists a function $V : \mathbb{R}^n \to \mathbb{R}_{\geq 0}$, functions $\alpha_1$, $\alpha_2$, $\sigma$ of class $\mathcal{K}_\infty$, and a continuous positive definite function $\rho$, such that (18) holds for all $x$, and so that $V(x(t, \xi, u))$ is Lebesgue measurable as a function of $t$, and (22) holds for all $\xi \in \mathbb{R}^n$, all input signals $u$, and all $t \geq 0$. Then the system is Integral Input-to-State Stable.

*Proof.* We apply Lemma 3.1 in [2] to the continuous positive definite function $\rho$, to conclude the existence of $\rho_1 \in \mathcal{K}_\infty$ and $\rho_2 \in \mathcal{L}$ such that:

$$\rho(r) \geq \rho_1(r)\,\rho_2(r) \tag{23}$$

holds for all $r \geq 0$. We define

$$\tilde{\rho}(r) := \rho_1(\alpha_2^{-1}(r))\,\rho_2(\alpha_1^{-1}(r))$$

and we pick a continuous positive definite and locally Lipschitz $\rho^*$ such that

$$\rho^*(r) \leq \tilde{\rho}_1(r/2)\,\tilde{\rho}_2(r)$$

for all $r$.

Thus, for all $x \in \mathbb{R}^n$ we have:

$$\rho(|x|) \geq \rho_1(|x|)\rho_2(|x|) \geq \rho_1(\alpha_2^{-1}(V(x)))\rho_2(\alpha_1^{-1}(V(x))) = \tilde{\rho}(V(x)). \tag{24}$$

From (22), taking into account (24), we have that

$$V(x(t, \xi, u)) \leq V(\xi) + \int_0^t \sigma(|u(s)|)\,ds - \int_0^t \tilde{\rho}(V(x(s, \xi, u)))\,ds \tag{25}$$

along all trajectories. We also note for future reference, that $V$ satisfies the following "no upper jumps" properties:

$$\limsup_{t \to t_0^+} V(x(t)) \leq V(x(t_0)) \tag{26}$$

$$\liminf_{t \to t_0^-} V(x(t)) \geq V(x(t_0)). \tag{27}$$

along all trajectories, for all $t_0 \geq 0$.

Now pick any initial state $\xi$ and any input $u$, and consider the (unique, since $\rho^*$ is locally Lipschitz) solution of the following initial value problem:

$$\dot{w} = \sigma(|u|) - \rho^*(w), \quad w(0) = V(\xi). \tag{28}$$

**Claim:**
$$V(x(t, \xi, u)) \leq w(t) \quad \forall t \geq 0. \tag{29}$$

In order to prove this, we first fix an arbitrary $\varepsilon > 0$, and we consider the initial value problem
$$\dot{w} = \sigma(|u|) - \rho^*(w), \quad w(0) = V(\xi) + \varepsilon. \tag{30}$$
We will we show that (29) holds for this modified problem, for each such $\varepsilon$. Then the result will follow for the original problem by letting $\varepsilon \to 0$ and using the fact that, for each $t$, the solution $w(t)$ depends continuously on $\varepsilon$. Note that, since $w(0) > 0$ for the modified problem, also $w(t) > 0$ for all $t$, by a comparison principle: $\dot{w}(t) \geq -\rho^*(w(t))$ for all $t$, and the equation $\dot{v} = -\rho^*(v)$ has unique solutions and has zero as an equilibrium.

Assume, by way of contradiction that there would exist some $t > 0$ such that $V(x(t,\xi,u)) \geq w(t)$. Let:
$$\tau := \inf\{t > 0 : V(x(t,\xi,u)) \geq w(t)\}. \tag{31}$$
To simplify notation, we write from now on $x(t)$ instead of $x(t,\xi,u)$. By definition of $\tau$ and (26), we have
$$w(\tau) \leq \limsup_{t \to \tau^+} V(x(t)) \leq V(x(\tau)). \tag{32}$$
Since $V(\xi) < w(0)$, necessarily $\tau > 0$.

Observe that $V(x(t)) < w(t)$ for all $t < \tau$, so in particular, using (27) and continuity of $w(t)$:
$$w(\tau) \geq \limsup_{t \to \tau^-} V(x(t)) \geq \liminf_{t \to \tau^-} V(x(t)) \geq V(x(\tau)). \tag{33}$$
Putting together (33) and (32) gives
$$w(\tau) = V(x(\tau)) = \lim_{t \to \tau^-} V(x(t)),$$
so there exists some $\delta > 0$ such that $\frac{w(t)}{2} < V(x(t)) < w(t)$ for every $t \in [\tau - \delta, \tau)$. We also have the fact that
$$\tilde{\rho}(r) \geq \tilde{\rho}_1(r)\tilde{\rho}_2(r) \geq \tilde{\rho}_1(s/2)\tilde{\rho}_2(s) = \rho^*(s)$$
whenever $s/2 \leq r \leq s$. We thus obtain a contradiction:
$$\begin{aligned}
V(x(\tau)) &\leq V(x(\tau-\delta)) + \int_{\tau-\delta}^{\tau} \sigma(|u(s)|)\,ds - \int_{\tau-\delta}^{\tau} \tilde{\rho}(V(x(s)))\,ds \\
&< w(\tau-\delta) + \int_{\tau-\delta}^{\tau} \sigma(|u(s)|)\,ds - \int_{\tau-\delta}^{\tau} \rho^*(w(s))\,ds \\
&= w(\tau).
\end{aligned}$$

This completes the proof of the claim.

By Corollary 4.3 in [2], associated to the positive definite continuous function $\rho^*$ there is some function $\beta$ of class $\mathcal{KL}$ with the following property: if $u$ is an input and $w$ is the solution of $\dot{w} = \sigma(|u|) - \rho^*(w)$ with initial condition $w_0$, then $w(t) \leq \beta(w_0, t) + \int_0^t \sigma(|u(s)|)\,ds$ for all $t$. It follows from our claim that
$$\alpha_1(|x(t,\xi,u)|) \leq V(x(t,\xi,u)) \leq w(t) \leq \beta(V(\xi), t) + \int_0^t \sigma(|u(s)|)\,ds. \tag{34}$$

This proves that the system is iISS. ∎

## 2.2  1 ⇒ 3

We show next that 1 implies 3. We must show the estimate in (7). By virtue of the converse Lyapunov characterization of IISS cited earlier, there exists a smooth function $V : \mathbb{R}^n \to \mathbb{R}_{\geq 0}$ so that (18) holds for suitable $\alpha_1$ and $\alpha_2$, such that (9) holds. Integrating, we obtain that (22) holds along all solutions. In particular, $|x(t, \xi, u)| \leq \kappa_1(|\xi|) + \kappa_2(\int_0^t \sigma(|u(s)|)ds)$, where we have defined $\kappa_1(s) = \alpha_1^{-1} \circ 2\alpha_2(s)$ and $\kappa_2(s) = \alpha_1^{-1} \circ 2s$. We apply Lemma 3.1 in [2] to the continuous positive definite function $\rho$, to conclude the existence of $\rho_1 \in \mathcal{K}_\infty$ and $\rho_2 \in \mathcal{L}$ such that (23) holds for all $r \geq 0$. Then,

$$V(x(t,\xi,u)) \leq V(\xi) + \int_0^t \sigma(|u(s)|)\,ds - \left[\int_0^t \rho_1(|x(s,\xi,u)|)\,ds\right] \left[\rho_2\left(\kappa_1(|\xi|) + \kappa_2\left(\int_0^t \sigma(|u(s)|)\,ds\right)\right)\right]. \quad (35)$$

Define now the function $\gamma$ as

$$\gamma(r) = \frac{1}{\rho_2(r)} - \frac{1}{\rho_2(0)}. \quad (36)$$

Notice that $\gamma$ is of class $\mathcal{K}_\infty$. Moreover, we also to define the following class $\mathcal{K}$ functions:

$$\begin{aligned}\chi_1(r) &= \max\{\gamma \circ 2\kappa_1(r), \alpha_2(r)\} \\ \chi_2(r) &= \max\{\gamma \circ 2\kappa_2(r), r\}. \end{aligned} \quad (37)$$

It follows from (35) and (37) that

$$\int_0^t \rho_1(|x(s,\xi,u)|)\,ds$$

$$\leq \left[\frac{1}{\rho_2(0)} + \gamma\left(\kappa_1(|\xi|) + \kappa_2\left(\int_0^t \sigma(|u(s)|)ds\right)\right)\right] \left[V(\xi) + \int_0^t \sigma(|u(s)|)\,ds\right]$$

$$\leq \left[\frac{1}{\rho_2(0)} + \gamma(2\kappa_1(|\xi|)) + \gamma \circ 2\kappa_2\left(\int_0^t \sigma(|u(s)|)ds\right)\right] \left[\alpha_2(|\xi|) + \int_0^t \sigma(|u(s)|)ds\right]$$

$$\leq \frac{\alpha_2(|\xi|) + \int_0^t \sigma(|u(s)|)ds}{\rho_2(0)} + \left(\chi_1(|\xi|) + \chi_2\left(\int_0^t \sigma(|u(s)|)\right)\right)^2$$

$$\leq \alpha_2(|\xi|)/\rho_2(0) + 2\chi_1(|\xi|)^2 + (\chi_2/\rho_2(0) + 2\chi_2^2) \circ \left(\int_0^t \sigma(|u(s)|)\right)$$

This is an estimate of type (6), as wanted.

## 2.3  3 ⇒ 2

As mentioned earlier, the fact that (7) together with forward completeness is a sufficient condition for O-GAS was already shown in [18]. So we must show that a UBEBS estimate holds.

We start by recalling a result in [3]. It is shown there that forward completeness of (1) implies (and is, in fact, equivalent to) there being an estimate of the following type along all trajectories:

$$|x(t,\xi,u)| \leq \kappa_1(t) + \kappa_2(|\xi|) + \kappa_3\left(\int_0^t \gamma(|u(s)|)\,ds\right) + c \quad (38)$$

holding for some $\kappa_1,\kappa_2,\kappa_3,\gamma$ of class $\mathcal{K}_\infty$ and some $c \geq 0$. We choose functions like this, and let $\alpha$, $\chi$, and $\sigma$ be as in the estimate (7). We also introduce
$$\delta := \max\{\gamma, \sigma\}.$$
We also define, for each $r \geq 0$:
$$m(r) := \sup\left\{|x(t,\xi,u)| : t \geq 0, \ |\xi| \leq r, \ \int_0^{+\infty} \delta(|u(s)|)\,ds \leq r\right\}.$$
Note that $m$ is a nondecreasing function. The main technical step is in showing that $m(r) < \infty$:

**Lemma 2.4** For each $r > 0$,
$$m(r) \leq \kappa_1\left(\frac{\chi(2r)}{\alpha(r)}\right) + \kappa_2(r) + \kappa_3(r) + c. \tag{39}$$

*Proof.* Pick any $r > 0$ and denote for simplicity the right-hand side of (39) as $M(r)$. Pick any state $\xi$ and input $u$ so that $|\xi| \leq r$ and $\int_0^{+\infty} \delta(|u(s)|)\,ds \leq r$. We need to show that, for all $t$, $|x(t,\xi,u)| \leq M(r)$. Assume that is not the case, so there is some $T > 0$ so that $|x(T,\xi,u)| > M(r)$. Let
$$\tau := \sup\{t \leq T : |x(t,\xi,u)| \leq r\},$$
so that $|x(t,\xi,u)| > r$ for all $t \in [\tau,T]$. It follows from (7) that
$$\begin{aligned}
\alpha(r)(T-\tau) &\leq \int_\tau^T \alpha(|x(s,\xi,u)|)\,ds \\
&\leq \int_0^T \alpha(|x(s,\xi,u)|)\,ds \\
&\leq \chi\left(r + \int_0^T \sigma(|u(s)|)\,ds\right) \\
&\leq \chi(2r).
\end{aligned}$$

With the notations
$$\tilde{T} := T - \tau \leq \frac{\chi(2r)}{\alpha(r)}, \quad \tilde{\xi} := x(\tau,\xi,u), \quad \tilde{u}(\cdot) := u(\cdot - \tau)$$
we have that
$$|x(T,\xi,u)| = |x(\tilde{T},\tilde{\xi},\tilde{u})| \leq \kappa_1(\tilde{T}) + \kappa_2(r) + \kappa_3(r) + c = M(r),$$
a contradiction. ∎

Now pick any $\xi$ and $u$, and let
$$r := \max\left\{|\xi|, \int_0^{+\infty} \delta(|u(s)|)\,ds\right\}.$$
By definition of $m$,
$$|x(t,\xi,u(\cdot))| \leq m(r) \leq \chi(r) + c$$
for all $t \geq 0$, where $\chi$ is any class-$\mathcal{K}_\infty$ function and $c$ is any constant such that $m(r) \leq \chi(r)+c$ for all $r$ (such $\chi$ and $c$ exist because $m$ is nodecreasing and finite-valued). With $\alpha := (2\chi)^{-1}$, we conclude that
$$\alpha(|x(t,\xi,u(\cdot))|) \leq r + \alpha(2c),$$
for all $t$, which gives us a UBEBS estimate, as wanted.

## 2.4  4 ⇒ 1

In order to prove the result, we need to appeal to the "small gain" argument used in the proof of the converse Lyapunov theorem for IISS in [2]. We review here the key technical steps needed from that argument, in a manner not stated explicitly in [2].

In general, for any system (1), and any given smooth $\mathcal{K}_\infty$ function $\varphi$, we consider the following auxiliary system:
$$\dot{x} = f(x, d\varphi(|x|)), \tag{40}$$
where we restrict the inputs $d$ to have values in the closed unit ball, i.e. $d(\cdot) : [0, +\infty) \to \bar{B}$, where $\bar{B}$ denotes the set $\{\mu \in \mathbb{R}^m : |\mu| \leq 1\}$. We let $\mathcal{M}_{\bar{B}}$ denote the set of such inputs. For each $\xi$ and $d$, we use $x^\varphi(t, \xi, d)$ to denote the trajectory of (40) corresponding to initial state $\xi$ and input $d$, defined in some maximal interval $[0, t_{\max}(\xi, d))$.

Suppose given, for the system (40), maps $\tilde{\alpha}$, $\gamma$, and $\varphi$ of class $\mathcal{K}_\infty$, such that $\gamma$ is smooth, and denote, for each $\xi \in \mathbb{R}^n$, each $t \geq 0$, and each $d \in \mathcal{M}_{\bar{B}}$,

$$z(t, \xi, d) := \tilde{\alpha}(|x^\varphi(t, \xi, d)|) - \int_0^t \gamma\Big(|d(s)|\,\varphi(|x^\varphi(s, \xi, d)|)\Big)\,ds$$

and for each $\xi$,
$$g(\xi) := \sup\{z(t, \xi, d) : t \in [0, t_{\max}(\xi, d)), d \in \mathcal{M}_{\bar{B}}\}$$

(possibly $= +\infty$). Suppose also given a $\beta_0 \in \mathcal{K}_\infty$ so that
$$z(t, \xi, d) \leq \beta_0(|\xi|) \quad \forall \xi, d, t \in [0, t_{\max}(\xi, d)). \tag{41}$$

Notice that, under this assumption, $g$ is finite-valued, and
$$\tilde{\alpha}(|\xi|) \leq g(\xi) \leq \beta_0(|\xi|)$$

for all $\xi \in \mathbb{R}^n$. Finally, assume that for each $0 < r_1 < r_2$ there is a $T(r_1, r_2) \geq 0$ such that the following property holds:

$$r_1 < |\xi| < r_2 \quad \Rightarrow \quad g(\xi) = \sup\{z(t, \xi, d) : 0 \leq t \leq T(r_1, r_2), d \in \mathcal{M}_{\bar{B}}\}. \tag{42}$$

Then, the proof of the main theorem in [2] contains a proof of this fact:

**Proposition 2.5** If the system (1) is O-GAS, and if the above assumptions hold, then the system (1) is IISS. □

We now apply this result to show that a system which satisfies an estimate (11) is necessarily IISS. Note that such a system is O-GAS. Assume without loss of generality that $\sigma = \gamma$ in (11), i.e., the system satisfies the following estimate:
$$\alpha(|x(t, \xi, u)|) \leq \beta(|\xi|, t) + \int_0^t \gamma(|u(s)|)\,ds + \gamma(\|u_{[0,t)}\|_\infty) \tag{43}$$

for all $t \geq 0$, all initial conditions $\xi \in \mathbb{R}^n$, and all measurable locally essentially bounded $u(\cdot)$. Let $\varphi$ be any smooth $\mathcal{K}_\infty$ function such that
$$\gamma(\varphi(s)) \leq \frac{\alpha(s)}{2} \quad \forall s \geq 0.$$

We will establish the properties needed for applying Proposition 2.5 with the same $\gamma$ and $\varphi$, $\tilde{\alpha} := \alpha/2$, and $\beta_0 := \beta(\cdot, 0)$. (The only minor technical problem in applying the result is the requirement that $\gamma$ be smooth; smoothness at the origin is needed due to the method of proof used in [2]. It is possible, however, to assume this fact with no loss of generality, employing the same trick as in [2] to produce a related system for which (11) holds with smooth $\gamma$.)

We start by establishing (41).

**Lemma 2.6** For all $\xi$, $d$, and $t \in [0, t_{\max}(\xi, d))$,

$$\alpha(|x^\varphi(t,\xi,d)|) \leq 2\beta_0(|\xi|) + 2\int_0^t \gamma\Big(|d(s)|\varphi(|x^\varphi(s,\xi,d)|)\Big)\,ds. \tag{44}$$

*Proof.* For all $\tau \in [0,t]$, write $u(\tau) := d(\tau)\varphi(|x^\varphi(\tau,\xi,d)|)$, so that

$$\gamma(|u(\tau)|) \leq \gamma(|\varphi(|x^\varphi(\tau,\xi,d)|)|) \leq \frac{\alpha(|x^\varphi(\tau,\xi,d)|)}{2} \tag{45}$$

for all $\tau$, and thus also

$$\gamma(\|u_{[0,\tau)}\|_\infty) \leq \frac{\alpha(\|x^\varphi_{[0,\tau)}\|_\infty)}{2}$$

for all $\tau$, where we are denoting by $\|x^\varphi_{[0,\tau)}\|_\infty$ the sup norm of $x^\varphi(\cdot,\xi,d)$ over the interval $[0,\tau]$. Now applying (43) with this $u$, we have

$$\alpha(|x^\varphi(\tau,\xi,d)|) \leq \beta(|\xi|,t) + \int_0^\tau \gamma(|u(s)|)\,ds + \frac{\alpha(\|x^\varphi_{[0,\tau)}\|_\infty)}{2}. \tag{46}$$

Now taking the supremum over $\tau \in [0,t]$ we obtain:

$$\alpha(\|x^\varphi_{[0,t)}\|_\infty) \leq \beta_0(|\xi|) + \int_0^t \gamma(|u(s)|)\,ds + \frac{\alpha(\|x^\varphi_{[0,t)}\|_\infty)}{2}$$

so that we conclude

$$\alpha(\|x^\varphi_{[0,t)}\|_\infty) \leq 2\beta_0(|\xi|) + 2\int_0^t \gamma(|u(s)|)\,ds.$$

Equation (44) follows from $\alpha(|x^\varphi(t,\xi,d)|) \leq \alpha(\|x^\varphi_{[0,t)}\|_\infty)$. ∎

In order to obtain the $T(r_1,r_2)$'s as in (42), we need a simple observation. For each $\xi$ and $d$, we let $G(\xi,d)$ denote the set of $t \in [0, t_{\max}(\xi,d))$ for which

$$|x^\varphi(t,\xi,d)| = \|x^\varphi_{[0,t)}\|_\infty$$

where we again use $\|x^\varphi_{[0,t)}\|_\infty$ to indicate the sup norm of $x^\varphi(\cdot,\xi,d)$ over the interval $[0,t]$.

**Lemma 2.7** For each $\xi$, $d$, and $t \in [0, t_{\max}(\xi,d))$, there is some $\tau \in G(\xi,d)$ such that $|x^\varphi(t,\xi,d)| = |x^\varphi(\tau,\xi,d)|$ and $\tau \leq t$.

*Proof.* Just let $\tau := \min\{s \in [0,t] : |x^\varphi(s,\xi,d)| \geq |x^\varphi(t,\xi,d)|\}$. By definition of $\tau$, $|x^\varphi(s,\xi,d)| < |x^\varphi(\tau,\xi,d)|$ for all $s \in [0,\tau)$. ∎

We use this remark as follows. Given any $\xi$ and $d$, and any $t$, we pick $\tau$ as in the Lemma. As $\tau \leq t$ and $\gamma$ is nonnegative,

$$-\int_0^t \gamma(|d(s)\varphi(|x^\varphi(s,\xi,d)|))\,ds \;\leq\; -\int_0^\tau \gamma(|d(s)\varphi(|x^\varphi(s,\xi,d)|))\,ds$$

so also $z(t,\xi,d) \leq z(\tau,\xi,d)$. Thus, for each $\xi$,

$$g(\xi) \;=\; \sup\{z(t,\xi,d) \,:\, t \in G(\xi,d),\, d \in \mathcal{M}_{\bar{B}}\}\,.$$

For $t \in G(\xi,d)$, the estimate (46) gives

$$\begin{aligned}\tilde{\alpha}(|x^\varphi(t,\xi,d)|) &= \alpha(|x^\varphi(t,\xi,d)|) - \frac{\alpha(\|x^\varphi_{[0,t)}\|_\infty)}{2} \\ &\leq \beta(|\xi|,t) + \int_0^t \gamma(|d(s)\varphi(|x^\varphi(s,\xi,d)|))\,ds\,,\end{aligned}$$

i.e., $z(t,\xi,d) \leq \beta(|\xi|,t)$. Since $g(\xi) \geq \tilde{\alpha}(|\xi|)$, we may pick as $T(r_1,r_2)$ any $t$ such that $\beta(r_2,t) < \tilde{\alpha}(r_1)$. This completes the verification of the conditions needed in order to apply Proposition 2.5.

## 3  Proof of Theorem 2

For technical reasons, it is convenient to introduce properties separating semiglobal behavior with respect to inputs or states respectively.

We say that a system is *semiglobal* IISS *with respect to inputs* if, for each $M > 0$ there are functions $\beta_M \in \mathcal{KL}$, and $\gamma_M$ and $\alpha_M$ in $\mathcal{K}_\infty$, such that the estimate (12) holds for all initial states $\xi$ and all those inputs $u$ for which $\|u\|_\infty \leq M$. The system is *semiglobal* IISS *with respect to initial states* if, for each $M > 0$ there are functions $\beta_M \in \mathcal{KL}$, and $\gamma_M$ and $\alpha_M$ in $\mathcal{K}_\infty$, such that the estimate (12) holds for all inputs $u$ and for all those initial states $\xi$ such that $|\xi| \leq M$. We could ask in these definitions "$\forall t \in [0, t_{\max}(\xi, u))$" or "$\forall t \geq 0$"; it amounts to the same thing since either property implies forward completeness, i.e. $t_{\max}(\xi, u) = +\infty$. Note also that if a system satisfies either of these properties then it is O-GAS (the argument is exactly the same as the one given earlier for semiglobal IISS).

The property of being semiglobally IISS with respect to inputs can be equivalently restated as the requirement that every saturated-input system

$$\dot{x} \;=\; f(x, \mathrm{sat}_M(u))\,, \tag{47}$$

where $\mathrm{sat}_M(u)$ indicates the projection of $u$ into the ball of radius $M$ in $\mathbb{R}^m$, be IISS, for each $M > 0$. Indeed, asking for each such system to be IISS amounts to asking that the required estimation functions exist.

We will first prove this:

**Proposition 3.1** The following facts are equivalent:

1. System (1) is IISS.

2. System (1) is semiglobally IISS with respect to inputs.

3. System (1) is semiglobally IISS with respect to initial states.

Then, we shall prove a variant of the equivalence of 4 and IISS in Theorem 1. We replace the estimate (11) by an estimate as follows:

$$\alpha(|x(t,\xi,u)|) \leq \beta_0(|\xi|) + \int_0^t \sigma(|u(s)|)\,ds + \gamma(\|u_{[0,t]}\|_\infty)\,, \tag{48}$$

understood as holding for some $\beta_0$, $\alpha$, $\sigma$, and $\gamma$ of class $\mathcal{K}_\infty$, along all trajectories. Of course, an estimate (11) implies also an estimate of this type, just using $\beta_0(r) := \beta(r,0)$. So an IISS system always satisfies this. The converse is less obvious:

**Proposition 3.2** System (1) is IISS if and only if it is 0-GAS and it satisfies an estimate (48).

The key calculation is contained in the following fact, to be proved after we show how the main conclusions follow from it.

**Lemma 3.3** Assume given three families of functions

$$\{\tilde{\beta}_M, M \in \mathbb{N}\} \subseteq \mathcal{KL}\,,\ \ \{\tilde{\sigma}_M, M \in \mathbb{N}\} \subseteq \mathcal{K}_\infty\,,\ \ \{\tilde{\gamma}_M, M \in \mathbb{N}\} \subseteq \mathcal{K}_\infty\,,$$

as well as two nondecreasing functions

$$\alpha_i : [0,\infty) \to [0,\infty)\,,\ \ i=1,2\,.$$

Then, there are a class-$\mathcal{KL}$ function $\beta$ and class-$\mathcal{K}_\infty$ functions $\gamma_1$, $\gamma_2$, $\delta_1$, and $\delta_2$ such that, for all $T \geq 0$,

$$\tilde{\beta}_{\lceil\alpha_1(R)+\alpha_2(S)\rceil}(R,T) + \tilde{\gamma}_{\lceil\alpha_1(R)+\alpha_2(S)\rceil}\left(\int_0^T \tilde{\sigma}_{\lceil\alpha_1(R)+\alpha_2(s)\rceil}(\varphi(s))\,ds\right)$$
$$\leq\ \beta(R,T) + \gamma_1(\alpha_1(R)) + \gamma_2(\alpha_2(S)) + \delta_1\left(\int_0^T \delta_2(\varphi(s))\,ds\right)$$

for all $R > 0$, all $S > 0$, and all measurable functions $\varphi : [0,T] \to \mathbb{R}_{\geq 0}$.

Here we are using $\lceil r \rceil$ to denote the "ceiling" of $r \in \mathbb{R}$, i.e. the smallest integer larger or equal than $r$. In the following we will also make use of the "floor" function, $\lfloor r \rfloor$, defined as the largest integer smaller or equal than $r$.

Let us see how Proposition 3.1 follows from this. Suppose first that the system is semiglobally IISS with respect to inputs. Thus there are families $\alpha_M$ and $\beta_M$ so that (12) holds whenever $\|u\|_\infty \leq M$, for all $\xi$ and $t$. Applying $\alpha_M^{-1}$ to both sides, we obtain:

$$|x(t,\xi,u)| \leq \tilde{\beta}_M(|\xi|,t) + \tilde{\gamma}_M\left(\int_0^t \tilde{\sigma}_M(|u(s)|)ds\right)\,, \tag{49}$$

where we defined $\tilde{\gamma}_M(r) = \alpha_M^{-1}(2r)$ and $\tilde{\beta}_M(r,t) = \alpha_M^{-1}(2\beta_M(r,t))$. We apply Lemma 3.3 to these families of maps (restricting attention to integer values of $M$), with $\alpha_1 \equiv 0$ and $\alpha_2(S) = S$. With the functions given by that Lemma, take now any $\xi$, $u$, and $t$. Let $S := \|u\|_\infty$, $M := \lceil S \rceil$, and $R := |\xi|$. As $\gamma_1(\alpha_1(R)) = \gamma_1(0) = 0$, we have

$$|x(t,\xi,u)| \leq \beta(|\xi|,t) + \gamma_2(\|u\|_\infty) + \delta_1\left(\int_0^t \delta_2(|u(s)|)\,ds\right)\,.$$

Applying to both sides $\alpha(\cdot)$, where $\alpha(3\delta_1(r)) = r$, we get rid of the $\delta_1$ function and recover the estimate in part 4 of Theorem 1, so that the system is indeed IISS.

Suppose now that, instead, the system is semiglobally IISS with respect to initial states. Thus there are families $\alpha_M$ and $\beta_M$ so that (12) holds whenever $|\xi| \leq M$, for all $u$ and $t$. Applying $\alpha_M^{-1}$ to both sides, we obtain again an estimate like (49). This time, we apply the Lemma with $\alpha_2 \equiv 0$ and $\alpha_1(R) = R$. With the functions given by the Lemma, take any $\xi$, $u$, and $t$. Let $S := \|u\|_\infty$, $R := |\xi|$, and $M := \lceil R \rceil$. As $\gamma_2(\alpha_1(S)) = 0$, we have

$$\begin{aligned}
|x(t,\xi,u)| &\leq \beta(|\xi|,t) + \gamma_1(|\xi|) + \delta_1\left(\int_0^t \delta_2(|u(s)|)\,ds\right) \\
&\leq \widehat{\gamma}(|\xi|) + \delta_1\left(\int_0^t \delta_2(|u(s)|)\,ds\right),
\end{aligned}$$

where $\widehat{\gamma}(r) = \beta(|r|,0) + \gamma_1(|r|)$. Applying to both sides $\alpha(\cdot)$, where $\alpha(2\delta_1(r)) = r$, we get rid of the $\delta_1$ function and recover an UBEBS-estimate; together with the fact that the system is 0-GAS, Theorem 1 guarantees that the system is indeed IISS. This ends the proof of Proposition 3.1.

Let us now prove Proposition 3.2. Suppose an estimate (48) holds for all trajectories. Because of Proposition 3.1, it is sufficient to show semiglobal IISS with respect to inputs, or equivalently that each system (47) is IISS. Fix any $M$. The estimate (48), applied to inputs bounded by $M$, gives us that the UBEBS property holds, with $c = \gamma(\|u_{[0,t]}\|_\infty)$ and $\gamma = \beta_0$. Thus the system (47) is IISS, as follows from Theorem 1 applied to that saturated system.

Finally, we prove Theorem 2. We assume that there are families $\alpha_M$ and $\beta_M$ so that (12) holds whenever $|\xi| \leq M$, $\|u\|_\infty \leq M$, and all $t$.' Applying $\alpha_M^{-1}$ to both sides, we obtain again an estimate like (49). We apply once more Lemma 3.3, now with $\alpha_1(R) = \alpha_2(R) = R$. With the functions given by the Lemma, take now any $\xi$, $u$, and $t$. Let $S := \|u\|_\infty$, $M := \lceil S \rceil$, and $R := |\xi|$. Then we obtain an estimate as follows along all solutions:

$$|x(t,\xi,u)| \leq [\beta(|\xi|,0) + \gamma_1(|\xi|)] + \gamma_2(\|u\|_\infty) + \delta_1\left(\int_0^t \delta_2(|u(s)|)\,ds\right).$$

Applying to both sides $\alpha(\cdot)$, where $\alpha(3\delta_1(r)) = r$, we once more eliminate the $\delta_1$ function and recover an estimate (48), so Proposition 3.2 gives us that the system is indeed IISS.

To summarize, we are only left with proving Lemma 3.3.

## 3.1 Proof of Lemma 3.3

Let us state several facts about $\mathcal{KL}$ and $\mathcal{K}_\infty$ functions, to be proved below.

**Lemma 3.4** Let $\{\tilde{\beta}_M\}_{M \in \mathbb{N}}$ be a family of class $\mathcal{KL}$ functions. Then there exist functions $\hat{\gamma}_1 \in \mathcal{K}_\infty$ and $\hat{\beta} \in \mathcal{KL}$ such that

$$\tilde{\beta}_M(r,t) \leq \hat{\gamma}_1(M)\,\hat{\beta}(r,t)$$

for all $r, t \geq 0$ and all $M > 0$.

**Lemma 3.5** Let $\{\tilde{\gamma}_M\}_{M \in \mathbb{N}}$ and $\{\tilde{\sigma}_M\}_{M \in \mathbb{N}}$ be two families of class $\mathcal{K}_\infty$ functions. Then, there exist three functions $\hat{\gamma}_2$, $\hat{\delta}_1$, and $\delta_2$ of class $\mathcal{K}_\infty$ such that, for all $T > 0$ and measurable functions $\varphi : [0,T] \to \mathbb{R}_{\geq 0}$,

$$\tilde{\gamma}_M\left(\int_0^T \tilde{\sigma}_M(\varphi(s))\,ds\right) \leq \hat{\gamma}_2(M)\,\hat{\delta}_1\left(\int_0^T \delta_2(\varphi(s))\,ds\right).$$

Lemma 3.3 follows immediately from these two lemmas, since $\lceil \alpha_1(R) + \alpha_2(S) \rceil \leq \alpha_1(R) + \alpha_2(S) + 1$ and using the facts that

$$\hat{\gamma}_1(\alpha_1(R) + \alpha_2(S) + 1) \leq \hat{\gamma}_1(3\alpha_1(R)) + \hat{\gamma}_1(3\alpha_2(S)) + \hat{\gamma}_1(3)$$

(and similarly for $\hat{\gamma}_2$) and $ab \leq a^2 + b^2$. One gets

$$\beta(R,T) = \hat{\gamma}_1(3)\,\hat{\beta}(R,T) + 2\,[\beta(R,T)]^2 \,,$$

$$\gamma_1(r) = [\hat{\gamma}_1(3r)]^2 + [\hat{\gamma}_2(3r)]^2 \,,$$

$$\delta_1(r) = \hat{\gamma}_2(3)\,\hat{\delta}_1(r) + 2\,[\delta_1(r)]^2 \,,$$

and the same $\delta_2$.

So we must prove these two lemmas. Let us collect first three useful facts about $\mathcal{KL}$ and $\mathcal{K}_\infty$ functions.

**Lemma 3.6** (Corollary 10 and Remark 11 in [18].) For each $\gamma \in \mathcal{K}_\infty$ there is a $\sigma \in \mathcal{K}_\infty$ such that $\gamma(rs) \leq \sigma(r)\sigma(s)$ for all $r, s \geq 0$. $\square$

**Lemma 3.7** (Proposition 7 in [18].) For each $\beta \in \mathcal{KL}$ there exist $\theta_1$ and $\theta_2$ in $\mathcal{K}_\infty$ such that $\beta(r,t) \leq \theta_1\left(\theta_2(s)e^{-t}\right)$ for all $r, t \geq 0$. $\square$

**Lemma 3.8** (Corollary 4.5 in [2].) Suppose that $\gamma : \mathbb{R}^2_{\geq 0} \to \mathbb{R}$ is such that $\gamma(\cdot, s) \in \mathcal{K}$ for each $s \in \mathbb{R}_{>0}$ and $\gamma(r, \cdot) \in \mathcal{K}$ for each $r \in \mathbb{R}_{>0}$. Then, there exists some function $\sigma \in \mathcal{K}$ such that

$$\gamma(r,s) \leq \sigma(r)\,\sigma(s)$$

for all $(x, y) \in (\mathbb{R}_{\geq 0})^2$. $\square$

We prove one more result of this type:

**Lemma 3.9** Let $\{\gamma_M\}_{M \in \mathbb{N}}$ be a family of class $\mathcal{K}_\infty$ functions. Then, there exists a $\sigma$ in $\mathcal{K}$ such that
$$\gamma_M(r) \leq \sigma(M)\,\sigma(r) \qquad \forall r \geq 0, \,\forall M \in \mathbb{N}\,.$$

*Proof.* We first make $\{\gamma_M\}$ into a monotonically increasing family of functions, by defining

$$\tilde{\gamma}_M(r) = \max_{i=1\ldots M} \gamma_i(r). \tag{50}$$

We then let $\hat{\gamma} : \mathbb{R}^2_{\geq 0} \to \mathbb{R}_{\geq 0}$ be as follows

$$\hat{\gamma}(s,r) = \begin{cases} \tilde{\gamma}_s(r) & \text{if } s \in \mathbb{N} \\ \tilde{\gamma}_{\lfloor s \rfloor}(r)(\lceil s \rceil - s) + \tilde{\gamma}_{\lceil s \rceil}(r)(s - \lfloor s \rfloor) & \text{if } s > 1 \text{ and } s \notin \mathbb{N} \\ \tilde{\gamma}_1(r)s & \text{if } s \in [0,1) \,. \end{cases} \tag{51}$$

By construction. $\hat{\gamma}$ is a function of class $\mathcal{KK}$, so we may apply to it Lemma 3.8, obtaining $\sigma$ so that

$$\gamma_M(r) \leq \tilde{\gamma}_M(r) \leq \hat{\gamma}(M,r) \leq \sigma(M)\,\sigma(r) \tag{52}$$

holds for all $M, r$. ∎

Let us prove Lemma 3.4. Let $\{\tilde{\beta}_M\}_{M \in \mathbb{N}}$ be a family of class $\mathcal{KL}$ functions. We apply Lemma 3.7 to each $\tilde{\beta}_M$, obtaining families $\{\theta_M^1\}$ and $\{\theta_M^2\}$ in $\mathcal{K}_\infty$ such that $\tilde{\beta}_M(r,t) \leq \theta_M^1 \left(\theta_M^2(r) e^{-t}\right)$ for all $r, t \geq 0$ and all $M$. Next we apply Lemma 3.9 to each of the families $\{\theta_M^1\}$ and $\{\theta_M^2\}$, obtaining functions $\sigma_1$ and $\sigma_2$, and to $\sigma_1$ we apply Lemma 3.6 to get $\theta \in \mathcal{K}_\infty$ such that:

$$\begin{aligned}
\tilde{\beta}_M(r,t) &\leq \sigma_1(M)\,\sigma_1\left(\sigma_2(M)\,\sigma_2(r)\,e^{-t}\right) \\
&\leq \underbrace{\sigma_1(M)\,\theta\left(\sigma_2(M)\right)}_{\tilde{\gamma}(M)}\,\underbrace{\theta\left(\sigma_2(r)\,e^{-t}\right)}_{\tilde{\beta}(s,t)}
\end{aligned}$$

as desired.

Finally, we prove Lemma 3.5. We apply Lemma 3.9 to each of the families $\{\tilde{\gamma}_M\}_{M \in \mathbb{N}}$ and $\{\tilde{\sigma}_M\}_{M \in \mathbb{N}}$, obtaining functions $\sigma_1$ and $\sigma_2$ and to $\sigma_1$ we apply Lemma 3.6 to get $\theta \in \mathcal{K}_\infty$ such that:

$$\begin{aligned}
\tilde{\gamma}_M\left(\int_0^T \tilde{\sigma}_M(\varphi(s))\,ds\right) &\leq \sigma_1(M)\,\sigma_1\left(\int_0^T \sigma_2(M)\,\sigma_2(\varphi(s))\,ds\right) \\
&\leq \sigma_1(M)\,\sigma_1\left(\sigma_2(M)\int_0^T \sigma_2(\varphi(s))\,ds\right) \\
&\leq \sigma_1(M)\,\theta\left(\sigma_2(M)\right)\,\theta\left(\int_0^T \sigma_2(\varphi(s))\,ds\right)
\end{aligned}$$

so the conclusions of Lemma 3.5 hold with with $\gamma_2(M) := \sigma_1(M)\,\theta\left(\sigma_2(M)\right)$, $\hat{\delta}_1 := \theta$, and $\delta_2 := \sigma_2$.

## 4 A Remark on Semiglobal ISS

For ISS, in contrast to iISS, the corresponding semiglobal property is strictly weaker, as shown by a counterexample in Section 4.1. If only inputs are restricted, however, a positive result exists, as we discuss next.

We say that a system is *semiglobally* ISS *with respect to inputs* if for each $M > 0$ there are functions $\beta_M \in \mathcal{KL}$, and $\gamma_M$ in $\mathcal{K}_\infty$, such that

$$|x(t,\xi,u)| \leq \beta_M(|\xi|,t) + \gamma_M(\|u\|_\infty) \tag{53}$$

holds for all initial states $\xi$ and all those inputs $u$ for which $\|u\|_\infty \leq M$.

We will establish the following result:

**Theorem 3** *System (1) is* ISS *if and only if it is semiglobally* ISS *with respect to inputs.*

*Proof.* We will prove the result exploiting the "LIM property" characterization of ISS given in [19]: ISS is equivalent to Lyapunov stability of the unforced system $\dot{x} = f(x,0)$ plus the following asymptotic gain condition on all solutions of (1):

$$\liminf_{t \to +\infty} |x(t,\xi,u)| \leq \gamma(\|u\|_\infty), \tag{54}$$

for some function $\gamma$ of class $\mathcal{K}_\infty$. Stability is clearly implied by semiglobal ISS with respect to inputs; in fact the unforced system is clearly O-GAS. Suppose that, for all $M$, the functions $\beta_M \in \mathcal{KL}$ and $\gamma_M$ in $\mathcal{K}_\infty$ are as in the definition. We apply lemma 3.9 to the family $\{\gamma_M\}$ to obtain a $\sigma$ as there, and define
$$\tilde{\gamma}(r) := \gamma(r+1)\,\gamma(r),$$
which is clearly of class $\mathcal{K}_\infty$. We claim that (54) holds. Indeed, pick any $\xi$ and $u$, and let $M := \lceil \|u\|_\infty \rceil$. Then, for all $t$,
$$|x(t,\xi,u)| \leq \beta_M(|\xi|,t) + \sigma(M)\,\sigma(\|u\|_\infty) \leq \beta_M(|\xi|,t) + \tilde{\gamma}(\|u\|_\infty)$$
from which (54) follows. ∎

### 4.1 Counter-example

We now provide a system which is "semiglobally ISS with respect to initial states" in the obvious sense but which is not ISS. Consider the system:
$$\begin{aligned} \dot{x}_1 &= -x_1(1-\sin x_2), \\ \dot{x}_2 &= -x_2 + u\,. \end{aligned} \tag{55}$$

*Claim 1: The system is not ISS.* Suppose that the system is ISS. Then there exists some $\beta \in \mathcal{KL}$ and some $\gamma \in \mathcal{K}$ such that
$$|x(t,\xi,u)| \leq \beta(|\xi|,t) + \gamma(\|u\|_\infty)$$
for all $t \geq 0$, all initial states $\xi$ and all $u$. Consequently,
$$\limsup_{t\to\infty} |x(t,\xi,u)| \leq \gamma(\|u\|_\infty) \tag{56}$$
for all $\xi$. Consider the initial state $\xi = (\gamma(\pi/2)+1, \pi/2)$ and the input given by $u(t) \equiv \pi/2$. The corresponding trajectory is
$$x_1(t) = \gamma(\pi/2)+1, \qquad x_2(t) = \pi/2, \qquad \forall\, t \geq 0.$$
Property (56) fails to hold for this trajectory. This shows that the system is not ISS.

*Claim 2:* There is a $\mathcal{KL}$-function $\beta$ such that, for each $M > 0$, there is a $\mathcal{K}$-function $\gamma_M$ such that
$$|x(t,\xi,u)| \leq \beta(|\xi|,t) + \gamma_M(\|u\|_\infty), \qquad \forall\, t \geq 0,$$
for all $|\xi| \leq M$ and all $u$.

The trajectories of the system are given by:
$$\begin{aligned} x_1(t) &= \xi_1 e^{\int_0^t -(1-\sin x_2(s))\,ds}, \\ x_2(t) &= \xi_2 e^{-t} + e^{-t}\int_0^t e^s u(s)\,ds\,. \end{aligned}$$

For the $x_2$ part, one has the following estimation:
$$|x_2(t)| \leq |\xi_2|e^{-t} + \|u\|_\infty. \tag{57}$$

The estimation for the $x_1$ is a bit more complicated than the $x_2$ part. First of all, we note that, for all $s \geq 0$,

$$\begin{aligned}
|\sin x_2(s)| &\leq \left|\sin[\xi_2 e^{-s} + e^{-s}\int_0^s e^\sigma u(\sigma)\,d\sigma]\right| \\
&= \left|\sin[\xi_2 e^{-s} + e^{-s}\int_0^s e^\sigma u(\sigma)\,d\sigma]\right| \\
&= \left|\sin(\xi_2 e^{-s})\cos(e^{-s}\int_0^s e^\sigma u(\sigma)\,d\sigma) + \cos(\xi_2 e^{-s})\sin(e^{-s}\int_0^s e^\sigma u(\sigma)\,d\sigma)\right| \\
&\leq |\sin(\xi_2 e^{-s})| + \left|\sin(e^{-s}\int_0^s e^\sigma u(\sigma)\,d\sigma)\right| \\
&\leq |\xi_2| + \left|\sin(e^{-s}\int_0^s e^\sigma u(\sigma)\,d\sigma)\right|,
\end{aligned}$$

so

$$\begin{aligned}
|x_1(t)| &\leq |\xi_1|e^{-t}e^{\int_0^t |\sin x_2(s)|\,ds} \\
&\leq |\xi_1|e^{|\xi_2|}e^{-t}e^{\int_0^t |\sin(e^{-s}\int_0^s e^\sigma u(\sigma)\,d\sigma)|\,ds}
\end{aligned}$$

Observe that

$$\left|\sin\left(e^{-s}\int_0^s e^\sigma u(\sigma)\,d\sigma\right)\right| \leq \min\left\{e^{-s}\int_0^t e^\sigma |u(\sigma)|\,d\sigma,\ 1\right\}.$$

Hence, when $\|u\|_\infty \leq 1/2$, one has

$$\int_0^t |\sin(e^{-s}\int_0^s e^\sigma u(\sigma)\,d\sigma)|\,ds \leq \int_0^t e^{-s}\int_0^s e^\sigma \|u\|_\infty\,d\sigma\,ds \tag{58}$$

$$\leq \int_0^t e^{-s}\frac{(e^s-1)}{2}\,ds \leq t/2, \tag{59}$$

and when $\|u\|_\infty \geq 1/2$, one has

$$\int_0^t \left|\sin(e^{-s}\int_0^s e^\sigma u(\sigma)\,d\sigma)\right|\,ds \leq \int_0^t ds = t.$$

Consequently, when $\|u\| \leq 1/2$,

$$|x_1(t)| \leq |\xi_1|e^{|\xi_2|}e^{-t}e^{t/2} = |\xi_1|e^{|\xi_2|}e^{-t/2},$$

and when $\|u\| \geq 1/2$,

$$|x_1(t)| \leq |\xi_1|e^{|\xi_2|}e^{-t}e^t = |\xi_1|e^{|\xi_2|}.$$

Let $\tilde\beta(s,t) = se^s e^{-t/2}$, and, for each $M > 0$, let $\tilde\gamma_M$ be any $\mathcal{K}$ function such that $\tilde\gamma_M(r) \geq Me^M$ for all $r \geq 1/2$. Then the above shows that

$$|x_1(t)| \leq \tilde\beta(|\xi|,t) + \tilde\gamma_M(\|u\|_\infty), \qquad \forall\, t \geq 0,$$

for all $\xi$ and $u$. Combining this with (57), one sees that Claim 2 is true.

# References

[1] D. Angeli, "Intrinsic robustness of global asymptotic stability", *Systems & Control Letters*, to appear. Preprints available at `http://www.dsi.unifi.it/~angeli/`

[2] D. Angeli, E. Sontag and Y. Wang, "A Lyapunov characterization of Integral Input-to-State Stability", *IEEE Transactions on Automatic Control*, submitted

[3] D. Angeli and E. Sontag, "Forward Completeness, Unboundedness Observability, and their Lyapunov characterizations", *Systems & Control Letters*, to appear.

[4] P.D. Christofides and A.R. Teel, "Singular perturbations and input-to-state stability," *IEEE Trans. Automat. Control* 41(1996): pp. 1645-1650.

[5] L. Grune, E. Sontag and F. Wirth, "Asymptotic stability equals exponential stability and ISS equals finite energy gain - if you twist your eyes", *Systems & Control Letters*, to appear.

[6] Hu, X.M., "On state observers for nonlinear systems," *Systems & Control Letters* 17 (1991), pp. 645–473.

[7] A. Isidori, "Global almost disturbance decoupling with stability for non minimum-phase single-input single-output nonlinear systems," *Systems & Control Letters* 28 (1996), pp. 115–122.

[8] Z. P. Jiang, A. Teel and L. Praly, "Small-gain theorem for ISS systems and applications," *Mathematics of Control, Signals, and Systems* 7, (1994), pp. 95-120.

[9] H.K. Khalil, *Nonlinear Systems*, Prentice-Hall, Upper Saddle River, NJ, second ed., 1996.

[10] M. Krstić and H. Deng, *Stabilization of Uncertain Nonlinear Systems*, Springer-Verlag, London, 1998.

[11] M. Krstić, I. Kanellakopoulos, and P. V. Kokotović, *Nonlinear and Adaptive Control Design*, John Wiley & Sons, New York, 1995.

[12] D. Liberzon, E. Sontag and Y. Wang, "On integral input-to-state stabilization", *Proc. American Control Conference*, June 1999, S. Diego, pp. 1598-1602

[13] Y. Lin, E. Sontag and Y. Wang, "A smooth converse Lyapunov theorem for robust stability", *SIAM Journal on Control and Optimization* 34, (1996), pp. 124-160.

[14] Hespanha, J.P, and A.S. Morse, "Certainty equivalence implies detectability," *Systems and Control Letters* 36(1999): pp. 1-13.

[15] R. Marino and P. Tomei, "Nonlinear output feedback tracking with almost disturbance decoupling," *IEEE Trans. Automat. Control* 44 (1999): pp. 18–28.

[16] L. Praly and Y. Wang, "Stabilization in spite of matched unmodelled dynamics and an equivalent definition of input-to-state stability," *Math. of Control, Signals, and Systems* 9 (1996): pp. 1-33.

[17] R. Sepulchre, M. Jankovic, and P.V. Kokotovic, "Integrator forwarding: a new recursive nonlinear robust design," *Automatica* 33 (1997): pp. 979-984.

[18] E. Sontag, "Comments on Integral variants of Input-to-State Stability", *Systems & Control Letters* 34, (1998), pp. 93-100

[19] E. Sontag and Y. Wang, "New characterizations of Input to State Stability", *IEEE Transactions on Automatic Control* 41, (1996), pp. 1283-1294.

[20] E. Sontag and Y. Wang, "On characterizations of the input-to-state stability property", *Systems & Control Letters* 24, (1995), pp. 351-359.

[21] E. Sontag and Y. Wang, "Notions of Input-to-Output Stability", submitted.

[22] Tsinias, J., "Input to state stability properties of nonlinear systems and applications to bounded feedback stabilization using saturation," *ESAIM Control Optim. Calc. Var.* 2 (1997): pp. 57-85.